\documentclass[11pt]{amsart}
 \usepackage{amssymb,amsfonts,amscd,amsbsy}
\usepackage[mathscr]{eucal}
\usepackage{url}

\newtheorem{thm}{Theorem}[section]
\newtheorem{lem}[thm]{Lemma}
\newtheorem{prop}[thm]{Proposition}
\newtheorem{cor}[thm]{Corollary}
\theoremstyle{definition}
\newtheorem{remark}[thm]{Remark}

\begin{document}

\title[Supercongruences for Ap{\'e}ry-like Numbers]
{Supercongruences for Ap{\'e}ry-like Numbers} 
 
\author{Robert Osburn and Brundaban Sahu}

\address{School of Mathematical Sciences, University College Dublin, Belfield, Dublin 4, Ireland}

\address{School of Mathematical Sciences, National Institute of Science Education and Research, Bhubaneswar 751005, India}

\email{robert.osburn@ucd.ie}

\email{brundaban.sahu@niser.ac.in}


\date{\today}

\begin{abstract}
It is known that the numbers which occur in Ap{\'e}ry's proof of the irrationality of $\zeta(2)$ have many interesting congruence properties while the associated generating function satisfies a second order differential equation. We prove supercongruences for a generalization of numbers which arise in Beukers' and Zagier's study of integral solutions of Ap{\'e}ry-like differential equations. 
\end{abstract}

\maketitle

\section{Introduction}

In the course of his work on proving the irrationality of $\zeta(2)$, Ap{\'e}ry introduced, for an integer $n \geq 0$, the following sequence of numbers \cite{apery}, \cite{vander}

$$\displaystyle B(n):=\sum_{j=0}^{n} \binom{n}{j}^2  \binom{n+j}{j}. $$

\noindent Several authors have subsequently investigated many interesting congruence properties for $B(n)$ and its generalizations. For example, Beukers \cite{beukers1} employed ``brute force methods" to prove

\begin{equation} \label{modp3r}
B(mp^{r} - 1) \equiv B(mp^{r-1} - 1) \pmod {p^{3r}}
\end{equation}

\noindent for any prime $p>3$ and integers $m$, $r \geq 1$. In \cite{sb}, Stienstra and Beukers related the $B(n)$'s to the $p$-th Fourier coefficient of a modular form. If we define

$$
 \eta^{6}(4z) =: \sum_{n=1}^{\infty} a(n) q^{n} 
$$

\noindent where

\begin{center}
$\displaystyle \eta(z) :=  q^{\frac{1}{24}} \prod_{n=1}^{\infty} (1 - q^{n})$
\end{center}

\noindent is the Dedekind eta function, $q:= e^{2{\pi}iz}$ and $z \in \mathbb{H}$, then they proved using the formal Brauer group of some elliptic K3-surfaces that for all odd primes $p$ and any $m$, $r \in \mathbb{N}$ with $m$ odd, we have

\begin{equation} \label{brauer}
B\Big(\frac{mp^{r} - 1}{2} \Big) - a(p) B\Big(\frac{mp^{r-1} - 1}{2}\Big) + (-1)^{\tfrac{p-1}{2}} p^2 B\Big(\frac{mp^{r-2} - 1}{2}\Big) \equiv  0 \pmod {p^r}.
\end{equation}

Congruence (\ref{modp3r}) is but one example of a general phenomena called {\it Supercongruences}. This term refers to the fact that congruences of this type are stronger than the ones suggested by formal group theory. It appeared in \cite{beukers1} and was the subject of the Ph.D. thesis of Coster \cite{cos}. In fact, Coster studied the generalized Ap{\'e}ry numbers (see, for example, Theorem 4.3.1 in \cite{cos})

\begin{equation*} 
u(n, A, B, \epsilon):= \sum_{j=0}^{n} \binom{n}{j}^{A} \binom{n+j}{j}^{B} \epsilon^{j}
\end{equation*} 

\noindent where $A$, $B \in \mathbb{N}$, $\epsilon=\pm 1$ and proved that 

\begin{equation} \label{gan1}
u(mp^r, A, B, \epsilon) \equiv u(mp^{r-1}, A, B, \epsilon) \pmod{p^{3r}}
\end{equation}

\noindent if $A \geq 3$ and

\begin{equation*} 
u(mp^r - 1, A, B, \epsilon) \equiv u(mp^{r-1} - 1, A, B, \epsilon) \pmod{p^{3r}}
\end{equation*}

\noindent if $B \geq 3$. Other examples of supercongruences have been observed in the context of number theory (see Chapter 11 in \cite{wom}), quantum theory \cite{k}, and algebraic geometry \cite{frv}. Currently, there is no systematic explanation for such congruences. Perhaps, as mentioned in \cite{sb}, they are related to formal Chow groups. 

It is known that the $B(n)$'s satisfy the recurrence relation

$$
(n+1)^2 B(n+1) = (11n^2 + 11n + 3) B(n) + n^2 B(n-1)
$$

\noindent for $n \geq 1$. This implies that the generating function 

$$\mathcal{B}(t) = \sum_{n=0}^{\infty} B(n) t^n$$

\noindent satisfies the differential equation

$$
L \mathcal{B}(t) = 0
$$

\noindent where 

\begin{center}
$\displaystyle L=t(t^2 + 11t - 1)\frac{d^2}{dt^2} + (3t^2 + 22t -1)\frac{d}{dt} + t-3$. 
\end{center}

\noindent In \cite{beukers3}, Beukers considers the differential equation

\begin{equation}\label{beukersequation}
((t^3+at^2+bt)F'(t))'+(t-\lambda)F(t)=0
\end{equation} 

\noindent where $a, b$ and $\lambda$ are rational parameters and asks for which values of these parameters
this equation has a solution in $\mathbb Z[[t]].$ This equation has a unique solution which is regular 
at the origin with $F(0)=0$ given by 

\begin{equation*}
F(t)=\sum^\infty_{n=0}u(n) t^n 
\end{equation*}
with $u(0)=1$ and satisfies the recurrence relation

\begin{equation*}
b(n+1)^2u(n+1)+(an^2+an-\lambda)u(n)+n^2u(n-1)=0
\end{equation*}

\noindent where $n \geq 1$. In \cite{zagier}, Zagier describes a search over a suitably chosen domain of $100$ million triples $(a, b, \lambda)$. He finds $36$ triples which yield an integral solution to (\ref{beukersequation}) and classifies seven as ``sporadic" cases (for a conjecture concerning the only cases where (\ref{beukersequation}) has an integral solution, see page 354 of \cite{zagier}). All seven cases (which include $B(n)$) have a binomial sum representation and a geometric origin. 

The purpose of this paper is to study congruences, akin to (\ref{gan1}), for a generalization of one of the ``sporadic" cases which can be expressed in terms of binomial sums and has a parametrization in terms of modular functions. For $A$, $B \in \mathbb{N}$, let

\begin{equation} \label{cab}
C(n, A, B):=\sum_{k=0}^{n} \binom {n} {k}^{A} \binom{2k}{k}^{B}.
\end{equation}

\noindent The first few terms in the sequence of positive integers $\{ C(n, 2, 1) \}_{n \geq 0}$ are as follows:

\begin{center}
$1$, $3$, $15$, $93$, $639$, $4653$, $35169$, $272835$, $\dotsc$
\end{center}

\noindent This sequence (see $A002893$ of Sloane \cite{sloane}) corresponds to a ``sporadic" case of Zagier (see $\# 8$ of Table 1 in \cite{zagier}) and has also appeared in the study of algebraic surfaces (see \cite{bea} or Part III of \cite{sb}), moments of Bessel functions arising in quantum field theory \cite{bbbg} and cooperative phenomena in crystals \cite{domb}. Our main result which is an analogue of (\ref{gan1}) is the following.

\begin{thm} \label{main} Let $A$, $B \in \mathbb{N}$ and $p > 3$ be a prime. For any $m$, $r \in \mathbb{N}$, we have

\begin{equation*} 
C(mp^r, A, B) \equiv C(mp^{r-1}, A, B) \pmod{p^{3r}}
\end{equation*}

\noindent if $A \geq 3$ and

\begin{equation*} 
C(mp^r, 2, B) \equiv C(mp^{r-1}, 2, B) \pmod{p^{2r}}.
\end{equation*}
\end{thm}

As a result of Theorem \ref{main}, we obtain a three-term congruence which is reminiscent of (\ref{brauer}). 

\begin{cor} \label{three} Let $A$, $B \in \mathbb{N}$ and $p > 3$ be a prime. Let $f$ be a non-cuspidal normalized Hecke eigenform of integer weight $k <4$ on $\Gamma_{0}(N)$ with character $\chi$ such that

\begin{equation} \label{nche}
f(z) =: \sum_{n=0}^{\infty} \gamma(n)q^n.
\end{equation}

\noindent Then for any $m$, $r \in \mathbb{N}$, we have

\begin{equation}\label{C}
C(mp^r, A, B) - \gamma(p) C(mp^{r-1}, A, B) + \chi(p) p^{k-1} C(mp^{r-2}, A, B) \equiv 0 \pmod {p^{3r + k - 4}}
\end{equation} 

\noindent if $A \geq 3$. If $k \geq 4$ and $A \geq 3$, this congruence is true modulo $p^{3r}$. If $A=2$ and $k<3$, it is true modulo $p^{2r+k-3}$. If $A=2$ and $k \geq 3$, it is true modulo $p^{2r}$. 

\end{cor}

The method of proof for Theorem \ref{main} is due to Coster (see page 50 of \cite{cos}). Namely, the idea is to rewrite the summands in (\ref{cab}) as products $g_{AB}(X,k)$ and $g_{AB}^{*}(X,k)$ (see Section 2), then exploit the combinatorial properties of these products. One then expresses (\ref{cab}) as two sums, one for which $p \mid k$ and the other for which $p \nmid k$. In the case $p \nmid k$, the sum vanishes modulo an appropriate power of $p$ while for $p \mid k$, the sum reduces to the required result. This approach can be used to prove supercongruences similar to Theorem \ref{main} in the remaining ``sporadic" cases. Finally, we would like to point out that Theorem 4.2 in \cite{ccs} follows from Theorem \ref{main} by taking $r=1$ and $B=0$ if $A \geq 3$ and from Lemma \ref{reduce} below once we use the identity

$$
C(n,2,0) = \binom{2n}{n}.
$$

The paper is organized as follows. In Section 2, we recall some properties of the products $g_{AB}(X,k)$ and $g_{AB}^{*}(X,k)$. In Section 3, we prove Theorem \ref{main} and Corollary \ref{three}. 

\section{Preliminaries}

We first recall the definition of two products and one sum and list some of their main properties. For more details, see Chapter 4 of \cite{cos}. For $A$, $B \in \mathbb{N}$ and integers $k$, $j \geq 1$ and $X$, we define

\begin{equation*} 
\displaystyle g_{AB}(X, k)=\prod_{i=1}^{k} \Biggl( 1 - \frac{X}{i} \Biggr)^{A} \Biggl(1 + \frac{X}{i} \Biggr)^{B},
\end{equation*}

\begin{equation*}
\displaystyle  g_{AB}^{*}(X, k)=\prod_{\substack{i=1 \\ p \nmid i}}^{k} \Biggl( 1 - \frac{X}{i} \Biggr)^{A} \Biggl(1 + \frac{X}{i} \Biggr)^{B},
\end{equation*}
 
\noindent and for a fixed prime $p>3$

\begin{equation*}
\displaystyle S_{j}(k)=\sum_{\substack{i=1 \\ p \nmid i}}^{k} \frac{1}{i^j}.
\end{equation*}

The following proposition provides some of the main properties of $g_{AB}(X, k)$, $g_{AB}^{*}(X,k)$ and $S_{j}(k)$. We note that parts (3) and (5) are straightforward to prove while (1), (2) and (4) require a short argument (see parts (i) and (ii) of Lemma 4.2.1 and parts (i), (ii) and (iv) of Lemma 4.2.5 in \cite{cos}).

\begin{prop} \label{prop}
For any $A$, $B$, $m \in \mathbb{N}$, $X \in \mathbb{Z}$ and integers $k, r \geq 1$, we have
\begin{enumerate}

\item $\displaystyle S_{j}(mp^r) \equiv 0 \pmod{p^r}$ for $j \not\equiv 0 \pmod {p-1}$,

\item $\displaystyle S_{2j-1}(mp^r) \equiv 0 \pmod{p^{2r}}$ for $j \not\equiv 0 \pmod{\frac{p-1}{2}}$,

\item $\displaystyle g_{AB}(pX, k) = g_{AB}^{*}(pX, k) g_{AB}(X, \bigl \lfloor \tfrac{k}{p} \bigr \rfloor)$,

\item $\displaystyle g_{AB}^{*}(X,k) \equiv 1 + (B-A) S_{1}(k)X + \tfrac{1}{2} \Bigl( (A-B)^2 S_{1}(k)^2 - (A+B) S_{2}(k) \Bigr) X^2  \pmod {X^3}$,

\item $\displaystyle \binom{n}{k}^{A} \binom{n+k}{k}^{B} =(-1)^{Ak} \Bigl ( \frac{n}{n-k} \Bigr)^{A} g_{AB}(n,k)$.

\end{enumerate}
\end{prop}

By (1), taking $j=1$ in (2) and (4) of Proposition \ref{prop}, we have 

\begin{equation} \label{key}
g_{AB}^{*}(mp^r, np^s) \equiv 1 \pmod{p^{r+2s}}
\end{equation} 

\noindent for any non-negative integers $m$, $n$, $r$ and $s$ with $s \leq r$. We now require a reduction result for one of the binomial coefficients occurring in $C(n, A, B)$.

\begin{lem} \label{reduce}

For a prime $p>3$ and integers $m \geq 0$, $r \geq 1$, we have

$$
\binom{2mp^r}{mp^r} \equiv \binom{2mp^{r-1}}{mp^{r-1}} \pmod{p^{3r}}.
$$

\end{lem}

\begin{proof} If $p \mid k$, then 

$$
\begin{aligned}
\binom{mp^r}{k} &= \binom{mp^{r-1}}{\tfrac{k}{p}} \prod_{\substack{\lambda=1 \\ p \nmid \lambda}}^{k} \Bigl( \frac{mp^r - \lambda}{\lambda} \Bigr) \\
& = \binom{mp^{r-1}}{\tfrac{k}{p}} (-1)^{k - \bigl \lfloor \tfrac{k}{p} \bigr \rfloor} \prod_{\substack{\lambda=1 \\ p \nmid \lambda}}^{k} \Bigl( 1 - \frac{mp^r}{\lambda} \Bigr) \\
& = \binom{mp^{r-1}}{\tfrac{k}{p}} g_{10}^{*} (mp^r, k).
\end{aligned}
$$

\noindent This implies 

$$
\binom{2mp^r}{mp^r} = \binom{2mp^{r-1}}{mp^{r-1}} g_{10}^{*}(2mp^r, mp^r). 
$$

\noindent By (\ref{key}), we have $g_{10}^{*}(2mp^r, mp^r) \equiv 1 \pmod {p^{3r}}$ and the result follows.

\end{proof}
 
\section{Proofs of Theorem \ref{main} and Corollary \ref{three}}

We are now in a position to prove Theorem \ref{main} and Corollary \ref{three}. \\

\noindent {\it{Proof of Theorem \ref{main}.}}
For integers $m$, $n$, $r \geq 1$ and $s \geq 0$ with $s \leq r$, we have

\begin{equation} \label{ord}
ord_{p} \binom{mp^r}{np^s}^{A} \geq A(r-s).
\end{equation}

\noindent Furthermore, by (3) and (5) of Proposition \ref{prop}, we have, for $s \geq 1$,

$$
\begin{aligned}
\binom{mp^r}{np^s}^{A} \binom{2np^s}{np^s}^{B} & = (-1)^{Anp^s} \Bigl ( \frac{mp^r}{mp^r-np^s} \Bigr)^{A} g_{A0}(mp^{r-1}, np^{s-1}) g_{A0}^{*}(mp^r, np^s) \binom{2np^s}{np^s}^{B} \\
& = \binom{mp^{r-1}}{np^{s-1}}^{A} g_{A0}^{*}(mp^r, np^s) \binom{2np^{s}}{np^s}^{B}. 
\end{aligned}
$$

\noindent By (\ref{key}) and Lemma \ref{reduce}, we have

$$
g_{A0}^{*}(mp^r, np^s) \equiv 1 \pmod{p^{r+2s}}
$$

\noindent and 

$$
\binom{2np^s}{np^s} \equiv \binom{2np^{s-1}}{np^{s-1}} \pmod{p^{3s}}.
$$

\noindent Thus

$$
\binom{mp^r}{np^s}^{A} \binom{2np^s}{np^s}^{B} \equiv \binom{mp^{r-1}}{np^{s-1}}^{A} \binom{2np^{s-1}}{np^{s-1}}^{B} \pmod{p^{min(A(r-s) + r + 2s, A(r-s) + 3s)}}.
$$

\noindent As $r \geq s$, we have $min(A(r-s) + r + 2s, A(r-s) + 3s)=A(r-s) + 3s$. If $A \geq 3$, then

\begin{equation} \label{summands}
\binom{mp^r}{np^s}^{A} \binom{2np^s}{np^s}^{B} \equiv \binom{mp^{r-1}}{np^{s-1}}^{A} \binom{2np^{s-1}}{np^{s-1}}^{B} \pmod{p^{3r}}.
\end{equation}

\noindent If $A=2$, then $A(r-s) + 3s \geq 2r$ and so 

\begin{equation} \label{A2}
\binom{mp^r}{np^s}^{A} \binom{2np^s}{np^s}^{B} \equiv \binom{mp^{r-1}}{np^{s-1}}^{A} \binom{2np^{s-1}}{np^{s-1}}^{B} \pmod{p^{2r}}.
\end{equation}

\noindent We now split $C(mp^r, A, B)$ into two sums, namely

$$
C(mp^r, A, B) = \sum_{\substack{k=0 \\ p \nmid k}}^{mp^r} \binom{mp^r}{k}^{A} \binom{2k}{k}^{B} + \sum_{\substack{k=0 \\ p \mid k}}^{mp^r} \binom{mp^r}{k}^{A} \binom{2k}{k}^{B}. 
$$

\noindent If $A \geq 3$, then the first sum vanishes modulo $p^{3r}$ using (\ref{ord}) and the result follows by (\ref{summands}). A similar argument is true for $A=2$ via (\ref{ord}) and (\ref{A2}).

\qed

\noindent {\it{Proof of Corollary \ref{three}.}} If $r \geq 2$, we have by Theorem \ref{main}

$$
C(mp^r, A, B) \equiv C(mp^{r-1}, A, B) \pmod{p^{3r}}
$$

\noindent if $A \geq 3$. Thus for $r \geq 3$

\begin{equation} \label{3rk4}
\chi(p) p^{k-1} C(mp^{r-1}, A, B) \equiv \chi(p) p^{k-1} C(mp^{r-2}, A, B) \pmod{p^{3r+k-4}}.
\end{equation}

\noindent The modular form (\ref{nche}) has the property that $\gamma(p)=1+ \chi(p) p^{k-1}$. This fact combined with (\ref{3rk4}) implies (\ref{C}) for $A \geq 3$ and $k<4$. Similarly, we have the result for the other cases.

\qed

\begin{remark}
We would like to mention an alternative approach to a weaker version of Corollary \ref{three} (with $p^{r}$ instead of $p^{2r}$) in the case $A=2$ and $B=1$ as this computation motivated Theorem \ref{main}. If we consider the modular function for $\Gamma_{0}(6)$

\begin{equation*} 
t(z)=\dfrac{\eta(z)^4\eta(6z)^8}{\eta(2z)^8 \eta(3z)^4},
\end{equation*}

\noindent then (see Theorems 5.1, 5.2 and 5.3 in \cite{cooper}, Table 2 in \cite{verrill}, Theorem 1 in \cite{yang} or Table 3, Case C in \cite{zagier})

\begin{equation*}
f(t):=\sum^\infty_{n=0}C(n, 2, 1)t^n =\dfrac{\eta(2z)^6\eta(3z)}{\eta(z)^3 \eta(6z)^2}.
\end{equation*}

\noindent If $\displaystyle \sigma_{\chi}(n):=\sum_{d \mid n} \chi(d) d^2$, then 

$$E(z)=-\frac{1}{9} + \sum_{n=1}^{\infty} \sigma_{\chi_{-3}}(n) q^n$$ 

\noindent is an Eisenstein series of weight $3$ on $\Gamma_{0}(3)$ with character $\chi_{-3}$. One can check that 

$$E(z) + 8E(2z)=-f(q) \dfrac{q \frac{dt}{dq}} {t}$$

\noindent and thus apply Theorem 1.1 in \cite{verrill}.

\end{remark}

\section*{Acknowledgements}
The authors were partially funded by Science Foundation Ireland 08/RFP/MTH1081. The authors would like to thank Frits Beukers for providing a copy of \cite{cos}, Dermot McCarthy for his comments and the referee for their careful reading of our paper.

\end{document}